\documentclass[12pt]{article}
\usepackage{amsmath,amsfonts,amssymb,amsthm}
\begin{document}

\newcommand{\1}{{{\bf 1}}}
\newcommand{\id}{{\rm id}}
\newcommand{\Hom}{{\rm Hom}\,}
\newcommand{\End}{{\rm End}\,}
\newcommand{\Res}{{\rm Res}\,}
\newcommand{\Image}{{\rm Im}\,}
\newcommand{\Ind}{{\rm Ind}\,}
\newcommand{\Aut}{{\rm Aut}\,}
\newcommand{\Ker}{{\rm Ker}\,}
\newcommand{\gr}{{\rm gr}}
\newcommand{\Der}{{\rm Der}\,}

\newcommand{\Z}{\mathbb{Z}}
\newcommand{\Q}{\mathbb{Q}}
\newcommand{\C}{\mathbb{C}}
\newcommand{\N}{\mathbb{N}}

\newcommand{\g}{\mathfrak{g}}
\newcommand{\h}{\mathfrak{h}}
\newcommand{\wt}{{\rm wt}\;}
\newcommand{\CR}{\mathcal{R}}
\newcommand{\D}{\mathcal{D}}
\newcommand{\E}{\mathcal{E}}
\newcommand{\Lie}{\mathcal{L}}
\newcommand{\z}{\bf{z}}
\newcommand{\bflam}{\bf{\lambda}}

\def \<{\langle}
\def \>{\rangle}
\def \be{\begin{equation}\label}
\def \ee{\end{equation}}
\def \bex{\begin{exa}\label}
\def \eex{\end{exa}}
\def \bl{\begin{lem}\label}
\def \el{\end{lem}}
\def \bt{\begin{thm}\label}
\def \et{\end{thm}}
\def \bp{\begin{prop}\label}
\def \ep{\end{prop}}
\def \br{\begin{rem}\label}
\def \er{\end{rem}}
\def \bc{\begin{coro}\label}
\def \ec{\end{coro}}
\def \bd{\begin{de}\label}
\def \ed{\end{de}}

\newtheorem{thm}{Theorem}[section]
\newtheorem{prop}[thm]{Proposition}
\newtheorem{coro}[thm]{Corollary}
\newtheorem{conj}[thm]{Conjecture}
\newtheorem{exa}[thm]{Example}
\newtheorem{lem}[thm]{Lemma}
\newtheorem{rem}[thm]{Remark}
\newtheorem{de}[thm]{Definition}
\newtheorem{hy}[thm]{Hypothesis}
\makeatletter
\@addtoreset{equation}{section}
\def\theequation{\thesection.\arabic{equation}}
\makeatother
\makeatletter


\begin{center}{\Large \bf On certain categories of modules
for twisted affine Lie algebras}\\
\vspace{0.5cm} Yongcun Gao\footnote{Partially supported by a grant
from Communication University of China} and Jiayuan
Fu\footnote{Partially supported by the National Science Foundation
of China (10726057)} \\
School of Science, Communication University of China, 100024\\
\end{center}

\begin{abstract}
In this paper, using generating functions we study two categories
${\mathcal{E}}$ and ${\mathcal{C}}$ of modules for twisted affine
Lie algebras $\hat{\g}[\sigma]$, which were firstly introduced and
studied in [Li] for untwisted affine Lie algebras. We classify
integrable irreducible $\hat{\g}[\sigma]$-modules in categories
${\mathcal{E}}$ and ${\mathcal{C}}$, where ${\mathcal{E}}$ is proved
to contain the well known evaluation modules and ${\mathcal{C}}$ to
unify highest weight modules, evaluation modules and their tensor
product modules. We determine the isomorphism classes of those
irreducible modules.
\end{abstract}

\baselineskip=16pt

\section{Introduction}
For affine Lie algebras, a very important class of modules is the
class of highest weight modules in the well known category
${\mathcal{O}}$, where highest weight integrable modules (of
nonnegative integral levels) (cf. [K]) have been the main focus. In
[C], Chari proved that every irreducible integrable
$\tilde{\g}$-module of a positive level with finite-dimensional
weight subspaces must be a highest weight module and that every
irreducible integrable $\hat{\g}$-module of level zero with
finite-dimensional weight subspaces must be finite-dimensional.

Another important class of modules is so called ``evaluation
modules" (of level zero) associated with a finite number of
$\g$-modules and with the same number of nonzero complex numbers,
studied by Chari and Pressley in [CP2] (cf. [CP1], [CP3]).  It was
proved in [CP2] (together with [C]) that every finite-dimensional
irreducible (integrable) $\hat{\g}$-module is isomorphic to an
evaluation module.  Furthermore, Chari and Pressley in [CP2] studied
the first time the tensor product module of an integrable highest
weight $\hat{\g}$-module with a (finite-dimensional) evaluation
$\hat{\g}$-module associated with finite-dimensional irreducible
$\g$-modules and distinct nonzero complex numbers.  A result, proved
in [CP2], is that such a tensor product module is also irreducible.
These irreducible (integrable) $\hat{\g}$-modules are greatly
different from highest weight modules and finite-dimensional modules
in many aspects.  In this way, a new family of irreducible
integrable $\hat{\g}$-modules were constructed.

With this new family of irreducible integrable $\hat{\g}$-modules
having been constructed, important problems are to determine the
isomorphism classes and to find a canonical characterization using
internal structures, instead of presenting them as tensor product
modules, and then to unify the new family with the family of highest
weight irreducible integrable modules toward a classification of all
the irreducible integrable $\hat{\g}$-modules.

In the paper [Li], all the above mentioned problems have been
completely solved for untwisted affine Lie algebras by exploiting
generating functions and formal calculus, which have played an
essential role in the theory of vertex operator algebras.

In [Li] two category $\E$ and ${\mathcal{C}}$  of modules for
untwisted affine Lie algebras $\hat{\g}$ are defined and studied. A
key result is a factorization which states that every irreducible
representation of $\hat{\g}$ in the category ${\mathcal{C}}$ can be
factorized canonically as the product of two representations such
that the first representation defines a restricted module and the
second one defines a module in the category $\E$.

We generalize all the results of [Li] to twisted affine Lie algebras
$\hat{\g}[\sigma]$. We classify integrable irreducible
$\hat{\g}[\sigma]$-modules in categories ${\mathcal{E}}$ and
${\mathcal{C}}$, where ${\mathcal{E}}$ is proved to contain the well
known evaluation modules and ${\mathcal{C}}$ to unify highest weight
modules, evaluation modules and their tensor product modules. We
determine the isomorphism classes of those irreducible modules.

This paper is organized as follows: In Section 2, we review the
notions of restricted module and integrable module and we prove a
complete reducibility theorem about the restricted integrable
$\hat{\g}[\sigma]$-modules. We also collect and restate certain
results on modules for tensor product algebras of two associative
algebras in [Li] . In Section 3, we study the category $\E$, and in
Section 4, we study the category ${\mathcal{C}}$. In Section 5, we
classify the irreducible integrable modules in the category
${\mathcal{C}}$.

\section{Category ${\mathcal{R}}$ of restricted $\hat{\g}[\sigma]$-modules}

First let us fix some formal variable notations
(see [FLM], [FHL], [LL]).
Throughout this paper,
$t, x, x_{1}, x_{2},\dots$
are independent mutually commuting formal variables.
We shall typically use $z, z_{1},z_{2},\dots$ for complex numbers.
For a vector space $U$, $U[[x_{1}^{\pm 1},\dots,x_{n}^{\pm 1}]]$
denotes the space of all formal (possibly doubly infinite) series
in $x_{1},\dots,x_{n}$ with coefficients in $U$, $U((x_{1},\dots,x_{n}))$
denotes the space of all formal (lower truncated) Laurent series
in $x_{1},\dots,x_{n}$ with coefficients in $U$
and $U[[x_{1},\dots,x_{n}]]$ denotes the space of all formal
(nonnegative) powers series
in $x_{1},\dots,x_{n}$ with coefficients in $U$.

We shall use the traditional binomial expansion convention:
For $m\in \Z$,
\begin{eqnarray}
(x_{1}\pm x_{2})^{m}=\sum_{i\ge 0}\binom{m}{i}(\pm 1)^{i}x_{1}^{m-i}x_{2}^{i}
\in \C[x_{1},x_{1}^{-1}][[x_{2}]].
\end{eqnarray}

Recall from [FLM] the formal delta function
\begin{eqnarray}
\delta(x)=\sum_{n\in \Z}x^{n}\in \C[[x,x^{-1}]].
\end{eqnarray}
Its fundamental property is that
\begin{eqnarray}
f(x)\delta(x)=f(1)\delta(x)\;\;\;\mbox{ for }f(x)\in \C[x,x^{-1}].
\end{eqnarray}
For any nonzero complex number $z$,
\begin{eqnarray}
\delta\left(\frac{z}{x}\right)=\sum_{n\in \Z}z^{n}x^{-n}\in \C[[x,x^{-1}]]
\end{eqnarray}
and we have
\begin{eqnarray}
f(x)\delta\left(\frac{z}{x}\right)
=f(z)\delta\left(\frac{z}{x}\right)\;\;\;\mbox{ for }f(x)\in \C[x,x^{-1}].
\end{eqnarray}
In particular,
\begin{eqnarray}
(x-z)\delta\left(\frac{z}{x}\right)=0.
\end{eqnarray}

Let $\g$ be a simple finite-dimensional Lie algebra equipped with a
nondegenerate symmetric invariant bilinear form $\<\cdot,\cdot\>$ so
that the squared length of the longest roots is $2$, and let
$\sigma$ be an automorphism of $\g$ of order $N,\,
\varepsilon=\mbox{exp}^\frac{2\pi i}{N}$. Then we have
$$\g=\bigoplus_{i=0}^{N-1}\g_{i},\,\,\mbox{where}\,\,\g_{i}=\{a\in \g|\,\sigma(a)=\varepsilon^{i}a\}.$$

Let $\hat{\g}[\sigma]$ be the corresponding twisted affine Lie
algebra, i.e.,
\begin{eqnarray}
\hat{\g}[\sigma]=\bigoplus_{i=0}^{N-1}\g_{i}\otimes
t^{\frac{i}{N}}\C[t,t^{-1}]\oplus \C {\bf k}
\end{eqnarray}
with the defining commutator relations
\begin{eqnarray}\label{edef-affine-comp}
[a\otimes t^{m+\frac{i}{N}},b\otimes t^{n+\frac{j}{N}}]
=[a,b]\otimes
t^{m+n+\frac{i+j}{N}}+(m+\frac{i}{N})\<a,b\>\delta_{m+\frac{i}{N},-n-\frac{j}{N}}{\bf
k},\label{eaffine1}
\end{eqnarray}
for $a\in\g_{i},b\in\g_{j},\; m,n\in \Z$, and with ${\bf k}$ as a
nonzero central element. A $\hat{\g}[\sigma]$-module $W$ is said to
be of {\em level} $\ell$ in $\C$ if the central element ${\bf k}$
acts on $W$ as the scalar $\ell$.

For $a\in\g_{i}$, form the generating function
\begin{eqnarray}
a(x)=\sum_{n\in \Z}(a\otimes
t^{n+\frac{i}{N}})x^{-n-\frac{i}{N}-1}\in \hat{\g}[\sigma]
[[x^{\frac{1}{N}},x^{-{\frac{1}{N}}}]].
\end{eqnarray}
 In terms of
generating functions the defining relations (\ref{eaffine1}) exactly
amount to
\begin{eqnarray}\label{edef-commutator}
[a(x_{1}),b(x_{2})]
=[a,b](x_{2})x_{1}^{-1}\left(\frac{x_{2}}{x_{1}}\right)^{\frac{i}{N}}\delta\left(\frac{x_{2}}{x_{1}}\right)+\<a,b\>
\frac{\partial}{\partial
x_{2}}\left[x_{1}^{-1}\left(\frac{x_{2}}{x_{1}}\right)^{\frac{i}{N}}\delta\left(\frac{x_{2}}{x_{1}}\right)\right]{\bf
k}.
\end{eqnarray}

Following the tradition (cf. [FLM], [LL]), for $a\in\g_{i},\;n\in\Z$
we shall use $a(n)$ for the corresponding operator associated to
$a\otimes t^{n+\frac{i}{N}}$ on $\hat{\g}[\sigma]$-modules. We have
the category ${\mathcal{R}}$ of the so-called restricted modules for
the affine Lie algebra $\hat{\g}[\sigma]$. A
$\hat{\g}[\sigma]$-module $W$ is said to be {\em restricted} (cf.
[K1]) if for any $w\in W,$
\begin{eqnarray}\label{erestricted-module}
a(n)w=0\;\;\;\mbox{ for }n\;\;\mbox{sufficiently large}.
\end{eqnarray}

Notice that in terms of generating functions, the condition (\ref{erestricted-module})
amounts to that
\begin{eqnarray}
a(x)w\in W((x^{\frac{1}{N}}))\;\;\;\mbox{ for } w\in W.
\end{eqnarray}
That is, a $\hat{\g}[\sigma]$-module $W$ is restricted if and only
if
\begin{eqnarray}\label{erestricted-condition}
a(x)\in \Hom (W,W((x^{\frac{1}{N}})))\;\;\;\mbox{ for }a\in \g.
\end{eqnarray}

Let $U$ be a $\g$-module and let $\ell$ be any complex number. Let
${\bf k}$ act on $U$ as the scalar $\ell$ and let
$$\hat{\g}[\sigma]_{+}=\sum_{n\in\N,\,i=0,\cdots,N-1,\,n+\frac{i}{N}>0}\g_{i}\otimes t^{n+\frac{i}{N}}$$act
trivially, making $U$ a $(\g_{0}\oplus\hat{\g}[\sigma]_{+}\oplus
\C{\bf k})$-module. Form the following induced
$\hat{\g}[\sigma]$-module
\begin{eqnarray}\label{einduced-module}
M_{\hat{\g}}(\ell,U)=U(\,\hat{\g}[\sigma]\,)\otimes_{\g_{0}\oplus\hat{\g}[\sigma]_{+}\oplus
\C{\bf k}}U.
\end{eqnarray}
It is clear that $M_{\hat{\g}}(\ell,U)$ is a restricted
$\hat{\g}[\sigma]$-module. We have also that the category
${\mathcal{R}}$ contains all the highest weight modules.

\bl{ldlm-basis} There are homogeneous elements $\{
a^{1},\dots,a^{r}\}$ of $\g=\bigoplus_{i=0}^{N-1}\g_{i}$ such that
\begin{eqnarray}
\g=\mbox{span}\{a^{1},a^{2},\dots,a^{r}\},\;
[a^{k}(m),a^{k}(n)]=0\;\mbox{ for }1\le k\le r,\; m,n\in \Z
\end{eqnarray}
and such that for $1\le k\le r$ and for any $n\in \Z$, $a^{k}(n)$
acts locally nilpotently on all integrable
$\hat{\g}[\sigma]$-modules. \el

It is similar to untwisted case, we have the following result:

\bt{tdlm} Let $\g$ be a finite-dimensional simple Lie algebra
equipped with the normalized Killing form.  Every nonzero restricted
integrable $\hat{\g}[\sigma]$-module is a direct sum of
(irreducible) highest weight integrable modules.  In particular,
every irreducible integrable $\hat{\g}[\sigma]$-module $W$ is a
highest weight integrable module.  \et

\begin{proof}
In view of the complete reducibility theorem in [K1] we only need to
show that every nonzero restricted integrable
$\hat{\g}[\sigma]$-module $W$ contains a highest weight integrable
(irreducible) submodule.

{\bf Claim 1:} {\em There exists a nonzero $u\in W$ such that
$\hat{\g}[\sigma]_{+}u=0.$}

For any nonzero $u\in W$, since $W$ is restricted,
$\hat{\g}[\sigma]_{+}u$ is finite-dimensional. For any $u\in W$, we
define $d(u)=\dim \hat{\g}[\sigma]_{+}u$. We need to prove that
there is a $0\ne u\in W$ such that $d(u)=0$.

Suppose that $d(u)>0$ for any $0\ne u\in W$. Take $0\ne u\in W$ such
that $d(u)$ is minimal. By Lemma \ref{ldlm-basis}, There are $\{
a^{1},\dots,a^{r}\}$ of $\g$ such that $a^{k}(n)$ locally
nilpotently act on $W$ for $k=1,\dots,r,\; n\in \Z$. Let
$l+\frac{j}{N}$ be the positive number such that $\g_{j}(l)u\ne 0$
and $\g_{i}(n)u=0$ whenever $n+\frac{i}{N}>l+\frac{j}{N}$. By the
definition of $l$, $a^{k}(l)u\ne 0$ for some $a^{k}\in\g_{j}, 1\le
k\le r$.

Notice that $a^{k}(l)^{s}u=0$ for some nonnegative integer $s$. Let
$m$ be the nonnegative integer such that $a^{k}(l)^{m}u\ne 0$ and
$a^{k}(l)^{m+1}u=0$. Set $v=a_{k}(l)^{m}u$. We will obtain a
contradiction by showing that $d(v)<d(u)$. First we prove that if
$a(n)u=0$ for some $a\in \g_{i},\; n+\frac{i}{N}>0$, then $a(n)v=0$.
In the following we will show by induction on $m$ that
$a(n)a_{k}(l)^mu=0$ for any $a\in \g_{i}$ and $m\ge 0$. If $m=0$
this is immediate. Now assume that the result holds for $m.$ Since
$[a,a^{k}](l+n)u=0$ (from the definition of $l$) and $a(n)u=0,$ by
the induction assumption that $a(n)a^{k}(l)^mu=0$ we have:
\begin{eqnarray}
[a,a^{k}](l+n)a^{k}(l)^{m}u=0,\;a(n)a^{k}(l)^{m}u=0.
\end{eqnarray}
Thus
\begin{eqnarray}
& &a(n)a^{k}(l)^{m+1}u
=[a(n),a^{k}(l)]a^{k}(l)^{m}u+a^{k}(l)a(n)a^{k}(l)^{m}u\nonumber\\
& &\ \ \ \ =[a,a^{k}](l+n)a^{k}(l)^{m}u+a^{k}(l)a(n)a^{k}(l)^{m}u
\nonumber\\
& &\ \ \ \ =0,
\end{eqnarray}
as required. In particular, we see that $a(n)v=a(n)a^{k}(l)^{m}u=0$.
Therefore, $d(v)\le d(u)$. Since $a^{k}(l)v=0$ and $a^{k}(l)u\ne 0$,
we have $d(v)< d(u)$, a contradiction.

{\bf Claim 2:} $W$ contains an irreducible
highest weight integrable submodule. Set
\begin{eqnarray}
\Omega(W)=\{u\in W\;|\hat{\g}[\sigma]_{+}u=0\}.
\end{eqnarray}
Then $\Omega(W)$ is a $\g_{0}$-submodule of $W$ and it is nonzero by
Claim 1. Since $a^{k}(0)$ for $k=1,\dots,r$ act locally nilpotently
on $\Omega(W)$, it follows from the PBW theorem that for any $u\in
\Omega(W)$, $U(\g_{0})u$ is finite-dimensional, so that $U(\g_{0})u$
is a direct sum of finite-dimensional irreducible $\g_{0}$-modules.
Let $u\in \Omega(W)$ be a highest weight vector for $\g_{0}$. It is
clear that $u$ is a singular vector for $\hat{\g}[\sigma]$. It
follows from [K] that $u$ generates an irreducible
$\hat{\g}$-module.
\end{proof}

We recall from [Li] the following result:

\bl{ldensity} Let $A_{1}$ and $A_{2}$ be associative algebras (with
identity) and let $U_{1}$ and $U_{2}$ be irreducible modules for
$A_{1}$ and $A_{2}$, respectively. If either $\End
_{A_{1}}U_{1}=\C$, or $A_{1}$ is of countable dimension, then
$U_{1}\otimes U_{2}$ is an irreducible $A_{1}\otimes A_{2}$-module.
\el

\bl{ltensor-decomposition-2}
Let $A_{1}$ and $A_{2}$ be associative algebras (with identity)
and let $U$ be an irreducible $A_{1}\otimes A_{2}$-module. Suppose that
$U$ as an $A_{1}$-module has an irreducible submodule $U_{1}$
and assume that either $A_{1}$ is of countable dimension
or $\End _{A_{1}}U_{1}=\C$. Then
$U$ is isomorphic to an $A_{1}\otimes A_{2}$-module of the form
$U_{1}\otimes U_{2}$ as in Lemma \ref{ldensity}.
\el

\bl{ltensor-decomposition}
Let $A_{1}$ and $A_{2}$ be associative algebras (with identity)
and let $W$ be an $A_{1}\otimes A_{2}$-module.
Assume that $A_{1}$ is of countable dimension and assume that
$W$ is a completely reducible $A_{1}$-module and
a completely reducible $A_{2}$-module.
Then $W$ is isomorphic to a direct sum
of irreducible $A_{1}\otimes A_{2}$-modules of the form
$U\otimes V$ with $U$ an
irreducible $A_{1}$-module and $V$ an irreducible $A_{2}$-module.
\el

\section{Category $\E$ of $\hat{\g}[\sigma]$-modules}

In this section, we study the category ${\mathcal{E}}$ of
$\hat{\g}[\sigma]$-modules, which is shown to include the well known
evaluation modules (of level zero). We show that every irreducible
integrable $\hat{\g}[\sigma]$-module in the category $\E$ is
isomorphic to a finite-dimensional evaluation module.

\bd{dcategoryE} {\em For the twisted affine Lie algebra
$\hat{\g}[\sigma]$, the category ${\mathcal{E}}$ is defined to
consists of $\hat{\g}[\sigma]$-modules $W$ for which there exists a
nonzero polynomial $p(x)\in \C[x]$, depending on $W$, such that
\begin{eqnarray}
p(x)a(x)w=0\;\;\;\mbox{ for }a\in \g,\; w\in W.
\end{eqnarray}}
\ed

\bl{l-level-E} The central element $\,{\bf k}$ of
$\,\hat{\g}[\sigma]$ acts as zero on any $\hat{\g}[\sigma]$-module
in the category ${\mathcal{E}}$. \el

\begin{proof} Let $W$ be a $\hat{\g}[\sigma]$-module in the category
${\mathcal{E}}$ with a nonzero polynomial $p(x)$ such that
$p(x)a(x)=0$ on $W$ for $a\in \g$. If $p(x)$ is a (nonzero)
constant, we have $a(x)=0$ for all $a\in \g$, i.e., $a(n)=0$ for
$a\in \g,\; n\in \Z$. In view of the commutator relation
(\ref{edef-affine-comp}) we see that ${\bf k}$ must be zero on $W$.
Assume that $p(x)$ is of positive degree, so that $p'(x)\ne 0$. Pick
$a,b\in \g_{0}$ such that $\<a,b\>=1$. (Notice that
$\<\cdot,\cdot\>$ is nondegenerate on $\g_{0}$.) Using the
commutator relations (\ref{edef-commutator}) we get
\begin{eqnarray}
0=p(x_{1})p(x_{2})[a(x_{1}),b(x_{2})]={\bf
k}p(x_{1})p(x_{2})x_{1}^{-1}\frac{\partial}{\partial
x_{2}}\delta\left(\frac{x_{2}}{x_{1}}\right).
\end{eqnarray}

\begin{eqnarray}
\Res_{x_{2}}p(x_{1})p(x_{2})x_{1}^{-1}\frac{\partial}{\partial
x_{2}}\delta\left(\frac{x_{2}}{x_{1}}\right)
&=&-\Res_{x_{2}}p(x_{1})p'(x_{2})x_{1}^{-1}\delta\left(\frac{x_{2}}{x_{1}}\right)\nonumber\\
&=&-p(x_{1})p'(x_{1}),
\end{eqnarray}
we get ${\bf k}p(x_{1})p'(x_{1})=0$, which implies that ${\bf k}=0$ on $W$.
\end{proof}

Next, we give some examples of $\hat{\g}[\sigma]$-modules in $\E$.
Let $U$ be a $\g$-module and let $z$ be a nonzero complex number.
With a fixed $z^\frac{1}{N}$, define an action of $\hat{\g}[\sigma]$
on $U$ by
\begin{eqnarray}
a(n)\cdot u&=&\sum_{i=0}^{N-1}z^{n+\frac{i}{N}}(a_{i}u)\;\;
\mbox{ for } a=\sum_{i=0}^{N-1}a_{i}\in \g,\; n\in \Z,\\
{\bf k}\cdot U&=&0.
\end{eqnarray}
Then $U$ equipped with the defined action is a
$\hat{\g}[\sigma]$-module (of level zero), which is denoted by
$U(z)$. If $U$ is an irreducible $\g$-module, it is clear that
$U(z)$ is an irreducible $\hat{\g}[\sigma]$-module. More generally,
let $U_{1},\dots,U_{r}$ be $\g$-modules and let $z_{1},\dots,z_{r}$
be nonzero complex numbers. Then the tensor product
$\hat{\g}[\sigma]$-module $\otimes_{k=1}^{r}U_{k}(z_{k})$ is called
an {\em evaluation module}.

We now show that the evaluation module
$\otimes_{k=1}^{r}U_{k}(z_{k})$ is in the category $\E$. For $a\in
\g_{i},\; u_{k}\in U_{k}(z_{k})=U_{k}$, we have
\begin{eqnarray}\label{eevaluation-generating}
a(x)(u_{1}\otimes \cdots\otimes u_{r}) &=&\sum_{n\in
\Z}\sum_{k=1}^{r}z_{k}^{n+\frac{i}{N}}x^{-n-\frac{i}{N}-1}
(u_{1}\otimes \cdots \otimes au_{k}\otimes \cdots \otimes u_{r})\nonumber\\
&=&\sum_{k=1}^{r}x^{-1}\left(\frac{z_{k}}{x}\right)^{\frac{i}{N}}\delta\left(\frac{z_{k}}{x}\right)
(u_{1}\otimes \cdots \otimes au_{k}\otimes \cdots \otimes u_{r}).
\end{eqnarray}
Since $(x-z_{k})\delta\left(\frac{z_{k}}{x}\right)=0$ for
$k=1,\dots,r$, we get $(x-z_{1})\cdots (x-z_{r})a(x)(u_{1}\otimes
\cdots\otimes u_{r})=0$. Thus we have proved:

\bl{lchari1} Let $U_{1},\dots, U_{r}$ be $\g$-modules and let
$z_{1},\dots,z_{r}$ be nonzero complex numbers. Then on the tensor
product $\hat{\g}[\sigma]$-module $U_{1}(z_{1})\otimes \cdots
\otimes U_{r}(z_{r})$,
\begin{eqnarray}
(x-z_{1})\cdots (x-z_{r})a(x)=0\;\;\;\mbox{ for }a\in \g.
\end{eqnarray}
In particular, the evaluation $\hat{\g}[\sigma]$-module
$U_{1}(z_{1})\otimes \cdots \otimes U_{r}(z_{r})$ is in the category
$\E$. \el

\bp{plevel-0-general} Let $\hat{\g}[\sigma]$ be a twisted affine Lie
algebra. Then for any irreducible $\g$-modules $U_{1},\dots,U_{r}$
and for any distinct nonzero complex numbers $z_{1},\dots,z_{r}$,
the tensor product $\hat{\g}[\sigma]$-module
$U_{1}(z_{1})\otimes\cdots \otimes U_{r}(z_{r})$ is irreducible. \ep

\begin{proof} Notice that the universal enveloping algebra $U(\hat{\g}[\sigma])$
is of countable dimension. It follows from Lemma \ref{ldensity} (and
induction) that $U_{1}(z_{1})\otimes \cdots \otimes U_{r}(z_{r})$ is
an irreducible module for the product Lie algebra
$\hat{\g}[\sigma]\oplus \cdots \oplus \hat{\g}[\sigma]$ ($r$
copies). Denote by $\pi$ the representation homomorphism map. For
$1\le j\le r$, denote by $\psi_{j}$ the  $j$-th embedding of
$\hat{\g}[\sigma]$ into $\hat{\g}[\sigma]\oplus \cdots \oplus
\hat{\g}[\sigma]$ ($r$ copies) and denote by $\psi$ the diagonal map
from $\hat{\g}[\sigma]$ to $\hat{\g}[\sigma]\oplus \cdots \oplus
\hat{\g}[\sigma]$ ($r$ copies). Then $\psi =\psi_{1}+\cdots
+\psi_{r}.$ We also extend the linear maps $\psi$ and
$\psi_{1},\dots,\psi_{r}$ on $\hat{\g}[\sigma][[x,x^{-1}]]$
canonically.

For $1\le j\le r$, set $p_{j}(x)=\prod_{k\ne
j}(x-z_{k})/(z_{j}-z_{k})$. Then
\begin{eqnarray}
p_{j}(x)\delta\left(\frac{z_{k}}{x}\right)
=p_{j}(z_{k})\delta\left(\frac{z_{k}}{x}\right)
=\delta_{j,k}\delta\left(\frac{z_{k}}{x}\right)
\end{eqnarray}
for $j,k=1,\dots,r$. Using (\ref{eevaluation-generating}) We have
that on $U_{1}(z_{1})\otimes \cdots \otimes U_{r}(z_{r})$,
\begin{eqnarray}
p_{j}(x)\pi \psi_{k}(a(x))=\delta_{j,k}\pi \psi_{k}(a(x))
\;\;\;\mbox{ for }1\le j, k\le r,\; a\in \g.
\end{eqnarray}
Thus on $U_{1}(z_{1})\otimes \cdots \otimes U_{r}(z_{r})$,
\begin{eqnarray}
p_{j}(x)\pi\psi (a(x))=\pi \psi_{j}(a(x)) \;\;\;\mbox{ for }1\le
j\le r,\;a\in \g,
\end{eqnarray}
which implies that
\begin{eqnarray}
\pi \psi_{j}(\hat{\g}[\sigma])\subset \pi\psi
(\hat{\g}[\sigma])\;\;\;\mbox{ for }j=1,\dots,r.
\end{eqnarray}
{}From this we have
\begin{eqnarray}
\pi \psi (\hat{\g}[\sigma]) =\pi \psi_{1}(\hat{\g}[\sigma])+ \cdots
+ \pi \psi_{r}(\hat{\g}[\sigma]).
\end{eqnarray}
It follows that $U_{1}(z_{1})\otimes\cdots \otimes U_{r}(z_{r})$ is
an irreducible $\hat{\g}[\sigma]$-module.
\end{proof}

\bp{pisomorphism-E} Let $U_{1},\dots,U_{r},\; V_{1},\dots,V_{s}$ be
nontrivial irreducible $\g$-modules and let $z_{1},\dots,z_{r}$ and
$\xi_{1},\dots,\xi_{s}$ be two groups of distinct nonzero complex
numbers. Then the $\hat{\g}[\sigma]$-module $U_{1}(z_{1})\otimes
\cdots \otimes U_{r}(z_{r})$ is isomorphic to $V_{1}(\xi_{1})\otimes
\cdots \otimes V_{s}(\xi_{s})$ if and only if $r=s$, $z_{j}=\xi_{j}$
and $U_{j}\cong V_{j}$ up to a permutation. \ep

\begin{proof} We only need to prove the ``only if'' part.
Let $U$ be any $\hat{\g}[\sigma]$-module in category $\E$. There
exists a (unique nonzero) monic polynomial $p(x)$ of least degree
such that $p(x)a(x)U=0$ for $a\in \g$. Clearly, isomorphic
$\hat{\g}[\sigma]$-modules in category $\E$ have the same monic
polynomial. If $U=U_{1}(z_{1})\otimes \cdots \otimes U_{r}(z_{r})$,
we are going to show that $p(x)=(x-z_{1})\cdots (x-z_{r})$ is the
associated monic polynomial. First, by Lemma \ref{lchari1} we have
that $p(x)a(x)=0$ on $U_{1}(z_{1})\otimes \cdots \otimes
U_{r}(z_{r})$ for $a\in \g$. Let $q(x)$ be any polynomial such that
$q(x)a(x)=0$ on $U_{1}(z_{1})\otimes \cdots \otimes U_{r}(z_{r})$
for $a\in \g$. Set $p_{j}(x)=\prod_{k\ne j}(x-z_{k})/(z_{j}-z_{k})$
for $j=1,\dots,r$ as in the proof of Proposition
\ref{plevel-0-general}. For $a\in \g_{i},\; u_{j}\in U_{j}$ with
$j=1,\dots,r$, we have
\begin{eqnarray*}
0&=&q(x)p_{j}(x)a(x)(u_{1}\otimes \cdots\otimes u_{r})\nonumber\\
&=&q(x)x^{-1}\left(\frac{z_{j}}{x}\right)^{\frac{i}{N}}\delta\left(\frac{z_{j}}{x}\right)
(u_{1}\otimes \cdots \otimes au_{j}\otimes\cdots\otimes u_{r})\nonumber\\
&=&q(z_{j})x^{-1}\left(\frac{z_{j}}{x}\right)^{\frac{i}{N}}\delta\left(\frac{z_{j}}{x}\right)
(u_{1}\otimes \cdots \otimes au_{j}\otimes\cdots\otimes u_{r}).
\end{eqnarray*}
Since each $U_{j}$ is a nontrivial $\g$-module, we must have
$q(z_{j})=0$ for $j=1,\dots,r$. Thus $p(x)$ divides $q(x)$. This
proves that $p(x)$ is the associated monic polynomial.

Assume that $U_{1}(z_{1})\otimes \cdots\otimes U_{r}(z_{r})$ is
isomorphic to $V_{1}(\xi_{1})\otimes\cdots\otimes V_{s}(\xi_{s})$
with $\Phi$ a $\hat{\g}[\sigma]$-module isomorphism map.  Then the
two tensor product modules must have the same associated monic
polynomial.  That is, $(x-z_{1})\cdots (x-z_{r})=(x-\xi_{1})\cdots
(x-\xi_{s})$.  Thus $r=s$ and up to a permutation $z_{j}=\xi_{j}$
for $j=1,\dots,r$. Assume that $z_{j}=\xi_{j}$ for $j=1,\dots,r$.
For $1\le j\le r,\; a\in \g_{i}$ and for $u_{k}\in U_{k},\; v_{k}\in
V_{k}$ with $k=1,\dots,r$, we have
\begin{eqnarray}
& &p_{j}(x)a(x)(u_{1}\otimes\cdots\otimes u_{r})
=x^{-1}\left(\frac{z_{j}}{x}\right)^{\frac{i}{N}}\delta\left(\frac{z_{j}}{x}\right)
(u_{1}\otimes \cdots \otimes au_{j}\otimes\cdots\otimes u_{r}),\\
& &p_{j}(x)a(x)(v_{1}\otimes\cdots\otimes v_{r})
=x^{-1}\left(\frac{z_{j}}{x}\right)^{\frac{i}{N}}\delta\left(\frac{z_{j}}{x}\right)
(v_{1}\otimes \cdots \otimes av_{j}\otimes\cdots\otimes v_{r}).
\end{eqnarray}
Then
\begin{eqnarray}\label{eithfactor}
\Phi(u_{1}\otimes \cdots \otimes au_{j}\otimes\cdots\otimes u_{r})
&=&\Res_{x}x^{-1}\delta\left(\frac{z_{j}}{x}\right)
\Phi(u_{1}\otimes \cdots \otimes au_{j}\otimes\cdots\otimes u_{r})\nonumber\\
&=&\Res_{x}\left(\frac{z_{j}}{x}\right)^{-\frac{i}{N}}\Phi\left(p_{j}(x)a(x)(u_{1}\otimes\cdots\otimes
u_{r})\right)\nonumber\\
&=&\Res_{x}\left(\frac{z_{j}}{x}\right)^{-\frac{i}{N}}p_{j}(x)a(x)\Phi(u_{1}\otimes\cdots\otimes u_{r})\nonumber\\
&=&\sigma_{j}(a)\Phi (u_{1}\otimes\cdots\otimes u_{r}),
\end{eqnarray}
where for $a\in \g_{i},\; v_{1}\in V_{1},\dots,v_{r}\in V_{r}$,
$$ \sigma_{j}(a)(v_{1}\otimes\cdots\otimes v_{r})=
(v_{1}\otimes \cdots \otimes av_{j}\otimes\cdots\otimes v_{r}). $$
It follows that $U_{j}$ is isomorphic to $V_{j}$. For example,
consider $j=1$. Pick up nonzero vectors $u_{k}\in U_{k}$ for $2\le
k\le r$. Using these vectors we get an embedding
$$\Theta: U_{1}\rightarrow U_{1}\otimes \cdots \otimes U_{r};\;\;\;\;
u\mapsto u\otimes u_{2}\otimes\cdots \otimes u_{r}.$$ On the other
hand, for any linear functions $f_{k}\in V_{k}^{*}$ for
$k=2,\dots,r$, we have a linear map
$$\Psi_{f_{2},\dots,f_{r}}:
V_{1}\otimes \cdots \otimes V_{r}\rightarrow V_{1};\;\; v_{1}\otimes
v_{2}\otimes \cdots \otimes v_{r}\mapsto f_{2}(v_{2})\cdots
f_{r}(v_{r})v_{1}.$$ Since $0\ne \Phi \Theta (U_{1})\subset
V_{1}\otimes \cdots \otimes V_{r}$, there exist linear functions
$f_{k}\in V_{k}^{*}$ for $k=2,\dots,r$ such that
$$\Psi_{f_{2},\dots,f_{r}}\Phi \Theta (U_{1})\ne 0.$$
By (\ref{eithfactor}), $\Psi_{f_{2},\dots,f_{r}}\Phi \Theta$ is a
nonzero $\g$-homomorphism from $U_{1}$ to $V_{1}$. It is an
isomorphism because both $U_{1}$ and $V_{1}$ are irreducible.
\end{proof}

We next classify finite-dimensional irreducible
$\hat{\g}[\sigma]$-modules in category $\E$. For $a_{i}\in \g_{i}$,
we have
\begin{eqnarray}
a_{i}(x)=\sum_{n\in \Z}(a_{i}\otimes
t^{n+\frac{i}{N}})x^{-n-\frac{i}{N}-1} =a_{i}\otimes
x^{-1}\delta\left(\frac{t}{x}\right)\left(\frac{t}{x}\right)^{\frac{i}{N}}.
\end{eqnarray}
For $f(x)\in \C [x],\; m\in \Z,\; a_{i}\in \g_{i}$, we have
\begin{eqnarray}
x^{m+\frac{i}{N}}f(x)a_{i}(x)=a_{i}\otimes
x^{m+\frac{i}{N}}f(x)x^{-1}\delta\left(\frac{t}{x}\right)\left(\frac{t}{x}\right)^{\frac{i}{N}}
=a_{i}\otimes
t^{m+\frac{i}{N}}f(t)x^{-1}\delta\left(\frac{t}{x}\right),
\end{eqnarray}
so that
\begin{eqnarray}
\Res_{x}x^{m+\frac{i}{N}}f(x)a_{i}(x)=a_{i}\otimes
t^{m+\frac{i}{N}}f(t).
\end{eqnarray}
It follows immediately that for any $\hat{\g}[\sigma] $-module $W$,
$f(x)a_{i}(x)W=0$ if and only if $(a_{i}\otimes
t^{\frac{i}{N}}f(t)\C[t,t^{-1}])W=0$.

For a nonzero polynomial $p(x)$, we define a subcategory
${\mathcal{E}}_{p(x)}$ of ${\mathcal{E}}$, consisting of
$\hat{\g}[\sigma]$-modules $W$ such that
\begin{eqnarray}
p(x)a(x)w=0\;\;\;\mbox{ for }a\in \g,\; w\in W.
\end{eqnarray}
Noticing Lemma \ref{l-level-E}, then a $\hat{\g}[\sigma]$-module in
the category $\E_{p(x)}$ exactly amounts to a module for the Lie
algebra $\sum_{i=0}^{N-1}\g_{i}\otimes
t^{\frac{i}{N}}\C[t,t^{-1}]/\sum_{i=0}^{N-1}\g_{i}\otimes
t^{\frac{i}{N}}p(t)\C[t,t^{-1}]$.

\bl{lEfactorization} Let $p(x)=(x-z_{1})\cdots (x-z_{r})$ with
$z_{1},\dots,z_{r}$ distinct nonzero complex numbers and with $k\in
\N$. Then any finite-dimensional irreducible
$\hat{\g}[\sigma]$-module $W$ in the category $\E_{p(x)}$ is
isomorphic to a $\hat{\g}[\sigma]$-module $U_{1}(z_{1})\otimes
\cdots \otimes U_{r}(z_{r})$ for some finite-dimensional irreducible
$\g$-modules $U_{1},\dots,U_{r}$. \el

\begin{proof} We have
$$\sum_{i=0}^{N-1}\g_{i}\otimes
t^{\frac{i}{N}}\C[t,t^{-1}]/\sum_{i=0}^{N-1}\g_{i}\otimes
t^{\frac{i}{N}}p(t)\C[t,t^{-1}]$$\\
$$=\bigoplus_{l=1}^{r}\left(\sum_{i=0}^{N-1}\g_{i}\otimes
t^{\frac{i}{N}}\C[t,t^{-1}]/\sum_{i=0}^{N-1}\g_{i}\otimes
t^{\frac{i}{N}}(t-z_{l})\C[t,t^{-1}]\right).$$

For any nonzero complex number $z$ and with a fixed
$z^{\frac{1}{N}}$, we define a map
$$\varphi:\sum_{i=0}^{N-1}\g_{i}\otimes
t^{\frac{i}{N}}\C[t,t^{-1}]/\sum_{i=0}^{N-1}\g_{i}\otimes
t^{\frac{i}{N}}(t-z)\C[t,t^{-1}]\rightarrow\g,$$ by
$$\sum_{i=0}^{N-1}a_{i}\otimes
t^{\frac{i}{N}}p_{i}(t)+\sum_{i=0}^{N-1}\g_{i}\otimes
t^{\frac{i}{N}}(t-z)\C[t,t^{-1}]\mapsto
\sum_{i=0}^{N-1}z^{\frac{i}{N}}p_{i}(z)a_{i}.$$   Then $\varphi$ is
an isomorphism. We have that a $\sum_{i=0}^{N-1}\g_{i}\otimes
t^{\frac{i}{N}}\C[t,t^{-1}]/\sum_{i=0}^{N-1}\g_{i}\otimes
t^{\frac{i}{N}}(t-z)\C[t,t^{-1}]$-module exactly amounts to an
evaluation $\hat{\g}[\sigma]$-module $U(z)$. Set
$$A_{l}=U\left(\sum_{i=0}^{N-1}\g_{i}\otimes
t^{\frac{i}{N}}\C[t,t^{-1}]/\sum_{i=0}^{N-1}\g_{i}\otimes
t^{\frac{i}{N}}(t-z_{l})\C[t,t^{-1}]\right)$$ for $l=1,\dots,r$.
Since $W$ is finite-dimensional, $W$ viewed as an $A_{l}$-module
contains a (finite-dimensional) irreducible submodule for which the
Schur Lemma holds. It now follows from Lemma
\ref{ltensor-decomposition-2} (and induction).
\end{proof}

\bp{pE-irreducible-module} Any finite-dimensional irreducible
$\hat{\g}[\sigma]$-module $W$ in the category $\E$ is isomorphic to
a $\hat{\g}[\sigma]$-module $U_{1}(z_{1})\otimes \cdots \otimes
U_{r}(z_{r})$ for some finite-dimensional $\g$-modules
$U_{1},\dots,U_{r}$ and for some distinct nonzero complex numbers
$z_{1},\dots,z_{r}$. \ep

\begin{proof} In view of Lemma \ref{lEfactorization},
it suffices to prove that $W$ is in the category $\E_{p(x)}$
with $p(x)$ a nonzero polynomial whose nonzero
roots are multiplicity-free.
Let $p(x)$ be the monic polynomial with the least degree
such that $p(x)a(x)W=0$ for $a\in \g$.
Notice that for any formal series $A(x)$ with coefficients
in any vector space and for any integer $m$,
$x^{m}A(x)=0$ if and only if $A(x)=0$.
In view of this we have $p(0)\ne 0$.
Thus
\begin{eqnarray}
p(x)=(x-z_{1})^{k_{1}}\cdots (x-z_{r})^{k_{r}},
\end{eqnarray}
where $z_{1},\dots,z_{r}$ are distinct nonzero complex numbers
and $k_{1},\dots,k_{r}$ are positive integers.

Let $I$ be the annihilating ideal of $W$ in $\hat{\g}[\sigma]$. Then
$\sum_{i=0}^{N-1}\g_{i}\otimes
t^{\frac{i}{N}}p(t)\C[t,t^{-1}]\subset I$ and $W$ is an irreducible
faithful $\hat{\g}[\sigma]/I$-module. Therefore $\hat{\g}[\sigma]/I$
is reductive (where we using the fact that $W$ is
finite-dimensional). Set $f(x)=(x-z_{1})\cdots (x-z_{r})$ and let
$k$ be the largest one among $k_{1},\dots,k_{r}$. We see that $p(x)$
is a factor of $f(x)^{k}$. It follows that the quotient space
$(\sum_{i=0}^{N-1}\g_{i}\otimes t^{\frac{i}{N}}f(t)\C[t,t^{-1}])/I$
is a solvable ideal of $\hat{\g}/I$. With $\hat{\g}[\sigma]/I$ being
reductive, $(\sum_{i=0}^{N-1}\g_{i}\otimes
t^{\frac{i}{N}}f(t)\C[t,t^{-1}])/I$ must be in the center of
$\hat{\g}[\sigma]/I$. From this we have that
$\sum_{i=0}^{N-1}\g_{i}\otimes
t^{\frac{i}{N}}f(t)\C[t,t^{-1}])\subset I$. This proves that
$f(x)a(x)W=0$ for $a\in \g$. Consequently, $f(x)=p(x)$, that is,
$k_{1}=\cdots =k_{r}=1$.
\end{proof}

\bp{pevaluation-simple-modules}
The irreducible integrable $\hat{\g}[\sigma]$-modules in the
category ${\mathcal{E}}$ up to isomorphism are exactly those
evaluation modules $U_{1}(z_{1})\otimes \cdots \otimes U_{r}(z_{r})$
where $U_{j}$ are finite-dimensional irreducible $\g$-modules and
$z_{j}$ are distinct nonzero complex numbers. \ep

\begin{proof} In view of Proposition \ref{pE-irreducible-module}
it suffices to prove that every irreducible integrable
$\hat{\g}[\sigma]$-module $W$ in the category ${\mathcal{E}}$ is
finite-dimensional. Since $W$ is in the category $\E$, there is a
nonzero polynomial $p(x)$ such that $(a\otimes t^\frac{i}{N}
p(t)\C[t,t^{-1}])W=0$ for $a\in \g_{i}, i=0,1,\cdots,N-1$. Let $I$
be the annihilating ideal of $W$ in $\hat{\g}[\sigma]$. Then
$\hat{\g}[\sigma]/I$ is finite-dimensional. Recall from Lemma
\ref{ldlm-basis} that there are homogeneous linearly generating
elements $\{ a^{1},\dots,a^{r}\}$ of $\g$ such that for any $1\le
k\le r,\; n\in \Z$, $a^{k}(n)$ acts locally nilpotently on $W$. Let
$0\ne w\in W$. Since $W$ is irreducible, we have
$W=U(\hat{\g}[\sigma])w=U(\hat{\g}[\sigma]/I)w$. In view of the PBW
theorem (for $\hat{\g}[\sigma]/I$ using a basis consisting of the
cosets of finitely many $a^{k}(n)$'s) we have that $W$ is
finite-dimensional, completing the proof.
\end{proof}

\section{Category ${\mathcal{C}}$ of $\hat{\g}$-modules}
In this section, we study the category ${\mathcal{C}}$ of
$\hat{\g}[\sigma]$-modules, which naturally unifies restricted
modules, evaluation modules and their tensor product modules. We
prove that the tensor product module of an irreducible restricted
$\hat{\g}[\sigma]$-module (in $\CR$) with an irreducible
$\hat{\g}[\sigma]$-module in $\E$ is irreducible. We also prove that
two such tensor product $\hat{\g}[\sigma]$-modules are isomorphic if
and only if the corresponding factors are isomorphic.

\bd{dcategoryC} {\em For a twisted affine Lie algebra
$\hat{\g}[\sigma]$, the category ${\mathcal{C}}$ is defined to
consist of $\hat{\g}[\sigma]$-modules $W$ such that there exists a
nonzero polynomial $p(x)\in \C[x]$, depending on $W$, such that
\begin{eqnarray}
p(x)a(x)\in \Hom (W,W((x^\frac{1}{N})))\;\;\;\mbox{ for }a\in \g.
\end{eqnarray}}
\ed

Clearly, every restricted $\hat{\g}[\sigma]$-module (in the category
${\mathcal{R}}$) and every $\hat{\g}[\sigma]$-module in the category
$\E$ are in ${\mathcal{C}}$. (Recall (\ref{erestricted-condition})
and Lemma \ref{lchari1}.)  It is also clear that ${\mathcal{C}}$ is
closed under tensor product of $\hat{\g}[\sigma]$-modules.  Thus,
tensor products of restricted $\hat{\g}[\sigma]$-modules with
evaluation $\hat{\g}[\sigma]$-modules belong to ${\mathcal{C}}$.
More specifically, {}from Lemma \ref{lchari1} we immediately have:

\bl{lchariC} Let $U_{1},\dots, U_{r}$ be $\g$-modules and let
$z_{1},\dots,z_{r}$ be nonzero complex numbers. For any restricted
$\hat{\g}[\sigma]$-module $W$, we have
\begin{eqnarray}
(x-z_{1})\cdots (x-z_{r})a(x)\in \Hom
(M,M((x^{\frac{1}{N}})))\;\;\;\mbox{ for }a\in \g,
\end{eqnarray}
where $M$ denotes the tensor product $\hat{\g}[\sigma]$-module
$W\otimes U_{1}(z_{1})\otimes \cdots \otimes U_{r}(z_{r})$. \el

\bt{ttensor-product-module} Let $W$ be an irreducible restricted
$\hat{\g}[\sigma]$-module (in the category ${\mathcal{R}}$) and let
$U$ be an irreducible $\hat{\g}[\sigma]$-module in the category
$\E$. Then the tensor product module $W\otimes U$ is irreducible.
\et

\begin{proof} Let $M$ be any nonzero
submodule of the tensor product $\hat{\g}[\sigma]$-module $W\otimes
U$. We must prove that $M=W\otimes U$. Since $W$ and $U$ are
irreducible $\hat{\g}[\sigma]$-modules and $U(\hat{\g}[\sigma])$ is
of countable dimension, by Lemma \ref{ldensity} $W\otimes U$ is an
irreducible $\hat{\g}[\sigma]\oplus \hat{\g}[\sigma]$-module. Now,
it suffices to prove that $M$ is a $\hat{\g}[\sigma]\oplus
\hat{\g}[\sigma]$-submodule of $W\otimes U$ and furthermore it
suffices to prove that
\begin{eqnarray}\label{eanm-M}
(a(n)\otimes 1)M\subset M\;\;\;\mbox{ for }a\in \g,\; n\in \Z.
\end{eqnarray}
(Notice that $(1\otimes a(n))w=a(n)w-(a(n)\otimes 1)w$ for $w\in M$.)

By the definition of $\E$ there exists a nonzero polynomial
$p(x)$ such that $p(x)a(x)=0$ on $U$ for all $a\in \g$, so that
\begin{eqnarray}\label{etensor=first}
p(x)(a(x)\otimes 1+1\otimes a(x))=p(x)a(x)\otimes 1
\;\;\mbox{ on }\;W\otimes U.
\end{eqnarray}
With $M$ being a $\hat{\g}[\sigma]$-submodule of the tensor product
module and with $W$ being a restricted module we have
\begin{eqnarray}
p(x)(a(x)\otimes 1+1\otimes a(x))M\subset M[[x,x^{-1}]],\;\;\;\;
p(x)(a(x)\otimes 1)M\subset (W\otimes U)((x)).
\end{eqnarray}
{}From this, using (\ref{etensor=first}) we have
\begin{eqnarray}
p(x)(a(x)\otimes 1)M\subset M((x))\;\;\;\mbox{ for }a\in \g.
\end{eqnarray}
Let $f(x)$ be the formal Laurent series of rational function $1/p(x)$
at zero, so that $f(x)\in \C((x))$.
Then we have
$$a(x)\otimes 1=(f(x)p(x))(a(x)\otimes 1)=f(x)(p(x)(a(x)\otimes 1))$$
on $M$. Consequently,
\begin{eqnarray}
(a(x)\otimes 1)M\subset M((x))\;\;\;\mbox{ for }a\in \g.
\end{eqnarray}
This proves (\ref{eanm-M}), completing the proof.
\end{proof}

The following result tells us when two $\hat{\g}[\sigma]$-modules of
the form $W\otimes U$ obtained in Theorem
\ref{ttensor-product-module} are isomorphic.

\bt{tidentification-2} Let $W_{1},W_{2}$ be irreducible
$\hat{\g}[\sigma]$-modules in category $\CR$ and let $U_{1}$ and
$U_{2}$ be irreducible $\hat{\g}[\sigma]$-modules in category $\E$.
Then the tensor product $\hat{\g}[\sigma]$-modules $W_{1}\otimes
U_{1}$ and $W_{2}\otimes U_{2}$ are isomorphic  if and only if
$W_{1}$ and $U_{1}$ are isomorphic to $W_{2}$ and $U_{2}$,
respectively. \et

\begin{proof} We only need to prove the ``only if'' part.
Let $F$ be a $\hat{\g}[\sigma]$-module isomorphism from
$W_{1}\otimes U_{1}$ onto $W_{2}\otimes U_{2}$. We have
\begin{eqnarray}\label{e2.36}
F(a(x)\otimes 1+1\otimes a(x))v=(a(x)\otimes 1+1\otimes a(x))F(v)
\;\;\;\mbox{ for }a\in \g,\; v\in W_{1}\otimes U_{1}.
\end{eqnarray}
Let $p(x)$ be a nonzero polynomial such that
$$p(x)a(x)U_{1}=0\;\;\mbox{ and }\;\; p(x)a(x)U_{2}=0\;\;\;\mbox{ for }a\in \g.$$
Using this and (\ref{e2.36}) we get
\begin{eqnarray}
p(x)F((a(x)\otimes 1)v)=p(x)(a(x)\otimes 1)F(v)
\;\;\;\mbox{ for }a\in \g,\; v\in W_{1}\otimes U_{1}.
\end{eqnarray}
We have
\begin{eqnarray}\label{efav}
F((a(x)\otimes 1)v)=(a(x)\otimes 1)F(v)
\;\;\;\mbox{ for }a\in \g,\; v\in W_{1}\otimes U_{1}.
\end{eqnarray}
Let $\{u_{2}^{i}\;|\;i\in S\}$ be a basis of $U_{2}$. Then
$W_{2}\otimes U_{2}=\coprod_{i\in S}W_{2}\otimes \C u_{2}^{i}$.
Denote by $\phi_{i}$ the projection of $W_{2}\otimes U_{2}$
onto $W_{2}\otimes \C u_{2}^{i}$.
We have
$$\phi_{i}(a(x)\otimes 1)w=(a(x)\otimes 1)\phi_{i}(w)
\;\;\;\mbox{ for }a\in \g,\; w\in
W_{2}\otimes U_{2},$$
so that
$$\phi_{i}F((a(x)\otimes 1)v)
=(a(x)\otimes 1)\phi_{i}F(v)\;\;\;\mbox{ for }a\in \g,\; v\in
W_{1}\otimes U_{1}.$$ Let $0\ne u_{1}\in U_{1}$. There exists an
$i\in S$ such that $\phi_{i}F\ne 0$ on $W_{1}\otimes \C u_{1}$. We
see that the map $\phi_{i}F$ gives rise to a nonzero
$\hat{\g}[\sigma]$-module homomorphism {}from $W_{1}\;(=W_{1}\otimes
\C u_{1})$ onto $W_{2}\;(=W_{2}\otimes \C u_{2}^{i})$. Because
$W_{1}$ and $W_{2}$ are irreducible, this nonzero homomorphism is an
isomorphism. This proves that $W_{1}$ is isomorphic to $W_{2}$.

{}From (\ref{e2.36}) and (\ref{efav}) we have
\begin{eqnarray}
F((1\otimes a(x))v)=(1\otimes a(x))F(v)
\;\;\;\mbox{ for }a\in \g,\; v\in W_{1}\otimes U_{1}.
\end{eqnarray}
Then using the same strategy, we see that $U_{1}$ is isomorphic to $U_{2}$.
\end{proof}

Theorem \ref{ttensor-product-module} gives us a construction of
irreducible $\hat{\g}[\sigma]$-modules in $\E$.  Naturally one wants
to know whether irreducible $\hat{\g}[\sigma]$-modules of the form
$W\otimes U$ as in Theorem \ref{ttensor-product-module} exhaust the
irreducible $\hat{\g}[\sigma]$-modules in the category
${\mathcal{C}}$ up to isomorphism. In the next section we shall
prove that this is true if we restrict ourselves to integrable
modules for $\hat{\g}[\sigma]$.

\section{Classification of
irreducible integrable $\hat{\g}[\sigma]$-modules in the category
${\mathcal{C}}$}

In this section we continue to study irreducible modules in the
category ${\mathcal{C}}$ for affine Lie algebra $\hat{\g}[\sigma]$.
As our main result we show that for the twisted affine Lie algebra
$\hat{\g}[\sigma]$, every irreducible integrable
$\hat{\g}[\sigma]$-module in the category ${\mathcal{C}}$ is
isomorphic to a tensor product of a highest weight integrable module
with a finite-dimensional evaluation module.

\bd{dbarE} {\em Let $W$ be any vector space. We set (see [Li]):
\begin{eqnarray}
\E(W)=\Hom (W,W((x^{\frac{1}{N}}))).
\end{eqnarray}
Define $\bar{\E}(W)$ to be the subspace of $(\End
W)[[x^{\frac{1}{N}},x^{-{\frac{1}{N}}}]]$, consisting of formal
series $a(x)$ such that $p(x)a(x)\in \Hom (W,W((x^{\frac{1}{N}})))$
for some nonzero polynomial $p(x)\in\C[x,x^{-1}]$ and define
$\bar{\E}_{0}(W)$ to be the subspace of $\bar{\E}(W)$ consisting of
the formal series $a(x)$ such that $p(x)a(x)=0$ for some nonzero
polynomial $p(x)$.} \ed

\br{rmodification} {\em For $a(x)\in (\End
W)[[x^{\frac{1}{N}},x^{-\frac{1}{N}}]]$, if $x^{m}f(x)a(x)\in \Hom
(W,W((x^{\frac{1}{N}})))$ for some integer $m$ and for some
polynomial $f(x)$, then $f(x)a(x)\in \Hom (W,W((x^{\frac{1}{N}})))$.
In view of this, if we need, we may assume that $p(0)\ne 0$ for a
nonzero polynomial $p(x)$ in Definition \ref{dbarE}. } \er

Let $\C(x)$ denote the algebra of rational functions of $x$.
We define $\iota_{x;0}$
to be the linear map from $\C(x)$ to $\C((x))$ such that
for $f(x)\in \C(x)$, $\iota_{x;0}(f(x))$ is the formal Laurent series
of $f(x)$ at $0$.
Notice that both $\C(x)$ and $\C((x))$ are
(commutative) fields.
The linear map $\iota_{x;0}$ is a field embedding.
If $p(x)$ is a polynomial
with $p(0)\ne 0$, then $\iota_{x;0}(p(x))\in \C[[x]]$.

\bd{dtildemap} {\em For a vector space $W$, we define a linear map
$$\psi_{{\mathcal{R}}}: \bar{\E}(W)\rightarrow \E(W)\;(=\Hom (W,W((x^{\frac{1}{N}}))))$$
by
\begin{eqnarray}\label{etildemap}
\psi_{{\mathcal{R}}}(a(x))w=\iota_{x;0}(f(x)^{-1})(f(x)a(x)w)
\;\;\;\mbox{ for }a(x)\in \bar{\E}(W), \;w\in W,
\end{eqnarray}
where $f(x)$ is any nonzero polynomial such that $f(x)a(x)\in \Hom
(W,W((x^{\frac{1}{N}})))$.} \ed

First of all, the map $\psi_{{\mathcal{R}}}$ is well defined;
the expression on the right-hand side
of (\ref{etildemap}) makes sense (which is clear)
and does not depend on the choice of $f(x)$.
Indeed,  let $0\ne f,g\in \C[x]$ be such that
$$f(x)a(x),\;\; g(x)a(x)\in \Hom (W,W((x^{\frac{1}{N}}))).$$
Set $h(x)=f(x)g(x)$. Then $h(x)a(x)\in
\Hom(W,W((x^{\frac{1}{N}})))$. For $w\in W$, we have
\begin{eqnarray*}
\iota_{x;0}(h(x)^{-1})(h(x)a(x)w)
&=&\iota_{x;0}(h(x)^{-1})f(x)(g(x)a(x)w)\\
&=&\iota_{x;0}(g(x)^{-1})(g(x)a(x)w).
\end{eqnarray*}
Similarly, we have
\begin{eqnarray*}
\iota_{x;0}(h(x)^{-1})(h(x)a(x)w)=\iota_{x;0}(f(x)^{-1})(f(x)a(x)w).
\end{eqnarray*}

The following is an immediate consequence of (\ref{etildemap}) and
the associativity law:
\bl{lbasic-property} For $a(x)\in
\bar{\E}(W)$, we have
\begin{eqnarray}
f(x)\psi_{\mathcal{R}}(a(x))=f(x)a(x),
\end{eqnarray}
where $f(x)$ is any nonzero polynomial such that $f(x)a(x)\in \Hom
(W,W((x^{\frac{1}{N}})))$. \el

Furthermore we have the following result:

\bp{pdecomposition}
For any vector space $W$, we have
\begin{eqnarray}\label{edecomposition}
\bar{\E}(W)=\E(W)\oplus \bar{\E}_{0}(W).
\end{eqnarray}
Furthermore, the linear map $\psi_{\mathcal{R}}$ from
$\bar{\E}(W)$ to $\E(W)$,
defined in Definition \ref{dtildemap},
is the projection map of $\bar{\E}(W)$ onto $\E(W)$,
i.e.,
\begin{eqnarray}\label{eprojections}
\psi_{\mathcal{R}}|_{\E(W)}=1\;\;\mbox{ and }
\;\; \psi_{\mathcal{R}}|_{\bar{\E}_{0}(W)}=0.
\end{eqnarray}
\ep

\begin{proof}
Let $a(x)\in \E(W)=\Hom (W,W((x^{\frac{1}{N}})))$. In Definition
\ref{dtildemap} we can take $f(x)=1$, so that
$\psi_{\mathcal{R}}(a(x))w=a(x)w$ for $w\in W$. Thus
$\psi_{\mathcal{R}}(a(x))=a(x)$.

Now, let $a(x)\in \bar{\E}_{0}(W)$. By definition there is a nonzero
polynomial $p(x)$ such that $p(x)a(x)=0$ on $W$, so that
$p(x)a(x)\in \Hom (W,W((x^{\frac{1}{N}})))$. From definition we have
$$\psi_{\mathcal{R}}(a(x))w=\iota_{x;0}(p(x)^{-1})(p(x)a(x)w)=0
\;\;\;\mbox{ for }w\in W.$$
Thus $\psi_{\mathcal{R}}(a(x))=0$.
This proves the property (\ref{eprojections})
and it follows immediately that the sum $\E(W)+ \bar{\E}_{0}(W)$ is
a direct sum.

Let $a(x)\in \bar{\E}(W)$ and let $0\ne f(x)\in \C[x]$ be such that
$f(x)a(x)\in \Hom (W,W((x^{\frac{1}{N}})))$. In view of Lemma
\ref{lbasic-property} we have
$f(x)\psi_{\mathcal{R}}(a(x))=f(x)a(x)$. Then
$f(x)(a(x)-\psi_{\mathcal{R}}(a(x)))=0$, which implies that
$a(x)-\psi_{\mathcal{R}}(a(x))\in \bar{\E}_{0}(W)$. Thus, $a(x)\in
\E(W)\oplus \bar{\E}_{0}(W)$. This proves that $\bar{\E}(W)\subset
\E(W)\oplus \bar{\E}_{0}(W)$, from which we have
(\ref{edecomposition}).
\end{proof}

\bd{dprojection-2}
{\em For a vector space $W$, we denote by $\psi_{\E}$ the projection map
of $\bar{\E}(W)$ onto $\E_{0}(W)$
with respect to the decomposition (\ref{edecomposition}).
For $a(x)\in \bar{\E}(W)$ we set
\begin{eqnarray}
\tilde{a}(x)&=&\psi_{{\mathcal{R}}}(a(x)),\\
\check{a}(x)&=&\psi_{\E}(a(x))=a(x)-\psi_{{\mathcal{R}}}(a(x))
=a(x)-\tilde{a}(x).
\end{eqnarray}}
\ed

{}From Lemma \ref{lbasic-property} we have
\begin{eqnarray}
f(x)\tilde{a}(x)&=&f(x)a(x),\\
f(x)\check{a}(x)&=&0
\end{eqnarray}
for any nonzero $f(x)\in \C[x]$ such that $f(x)a(x)\in
\Hom(W,W((x^{\frac{1}{N}})))$.

The following result relates the actions of $\psi_{\mathcal{R}}(a(x))$ and
$a(x)$ on $W$:

\bl{lconnection}
For $a(x)\in \bar{\E}(W),\; n\in \Z,\;w\in W$, we have
\begin{eqnarray}
\psi_{{\mathcal{R}}}(a(x))(n)w=\sum_{i=0}^{r}\beta_{i}a(n+i)w
\end{eqnarray}
for some $r\in \N,\;\beta_{1},\dots,\beta_{r}\in \C$,
depending on $a(x),w$ and $n$, where
\begin{eqnarray}
\psi_{{\mathcal{R}}}(a(x))=\sum_{n\in \Z}\psi_{{\mathcal{R}}}(a(x))(n)x^{-n-1}.
\end{eqnarray}
\el

\begin{proof}
Let $p(x)$ be a polynomial with $p(0)\ne 0$ such that $p(x)a(x)\in
\Hom (W,W((x^{\frac{1}{N}})))$. Then $x^{k}p(x)a(x)w\in
W[[x^{\frac{1}{N}}]]$ for some nonnegative integer $k$. Assume that
\begin{eqnarray}
\iota_{x;0}(1/p(x))=\sum_{i\ge 0}\alpha_{i}x^{i}\in \C[[x]].
\end{eqnarray}
Noticing that $\Res_{x}x^{k+m}p(x)a(x)w=0$ for $m\ge 0$,
we have
\begin{eqnarray}
\psi_{{\mathcal{R}}}(a(x))(n)w&=&\Res_{x}x^{n}\psi_{{\mathcal{R}}}(a(x))w
=\Res_{x}x^{n}\iota_{x;0}(1/p(x))(p(x)a(x)w)\nonumber\\
&=&\Res_{x}\sum_{0\le i\le k-n-1}\alpha_{i}x^{n+i}(p(x)a(x)w)\nonumber\\
&=&\Res_{x}\left(\sum_{0\le i\le k-n-1}\alpha_{i}x^{n+i}p(x)a(x)\right)w.
\end{eqnarray}
Then it follows immediately.
\end{proof}

\bl{lcommutatorrelation}
Let $W$ be a vector space as before and let
$$a(x), b(x)\in \bar{\E}(W),\;
c_{0}(x), c_{1}(x)\in (\End W)[[x^{\frac{1}{N}},x^{-\frac{1}{N}}]]$$
be such that on $W \,\,\mbox{for some}\,\,i\in\N$,
\begin{eqnarray}\label{ecommutatorold}
[a(x_{1}),b(x_{2})]=c_{0}(x_{2})x_{1}^{-1}\left(\frac{x_{2}}{x_{1}}
\right)^{\frac{i}{N}}\delta\left(\frac{x_{2}}{x_{1}}\right)
+c_{1}(x_{2})\frac{\partial}{\partial x_{2}}\left[
x_{1}^{-1}\left(\frac{x_{2}}{x_{1}}\right)^{\frac{i}{N}}
\delta\left(\frac{x_{2}}{x_{1}}\right)\right].
\end{eqnarray}
Then $c_{0}(x), c_{1}(x)\in \bar{\E}(W)$ and
\begin{eqnarray}\label{ecommutatornew}
[\tilde{a}(x_{1}),\tilde{b}(x_{2})]
=\tilde{c}_{0}(x_{2})x_{1}^{-1}\left(\frac{x_{2}}{x_{1}}\right)^{\frac{i}{N}}\delta\left(\frac{x_{2}}{x_{1}}\right)
+\tilde{c}_{1}(x_{2})\frac{\partial}{\partial
x_{2}}\left[x_{1}^{-1}\left(\frac{x_{2}}{x_{1}}\right)^{\frac{i}{N}}\delta\left(\frac{x_{2}}{x_{1}}\right)\right].
\end{eqnarray}
\el

\begin{proof}
From (\ref{ecommutatorold}) we get
\begin{eqnarray}
&
&\Res_{x_{1}}\left(\frac{x_{2}}{x_{1}}\right)^{-\frac{i}{N}}[a(x_{1}),b(x_{2})]
=c_{0}(x_{2})+\frac{i}{N}x_{2}^{-1}c_{1}(x_{2}),\\
&
&\Res_{x_{1}}(x_{1}-x_{2})\left(\frac{x_{2}}{x_{1}}\right)^{-\frac{i}{N}}[a(x_{1}),b(x_{2})]
=c_{1}(x_{2}).
\end{eqnarray}
Then it is clear that $c_{k}(x)\in \bar{\E}(W)$ for $k=0,1$, since
$b(x)\in \bar{\E}(W)$.

Let $0\ne f(x)\in \C[x]$ be such that
$$f(x)a(x),\;f(x)b(x),\; f(x)c_{0}(x),f(x)c_{1}(x)\in \Hom (W,W((x))),$$
so that
$$f(x)a(x)=f(x)\tilde{a}(x),\; f(x)b(x)=f(x)\tilde{b}(x),\;
f(x)c_{k}(x)=f(x)\tilde{c}_{k}(x)$$ for $k=0,1$. Then multiplying
both sides of (\ref{ecommutatorold}) by $f(x_{1})f(x_{2})$ we obtain
\begin{eqnarray}
&&f(x_{1})f(x_{2})[\tilde{a}(x_{1}),\tilde{b}(x_{2})]\\
&&=\left[f(x_{1})f(x_{2})\tilde{c}_{0}(x_{2})
+f(x_{1})f(x_{2})\tilde{c}_{1}(x_{2})\frac{\partial}{\partial
x_{2}}\right]x_{1}^{-1}\left(\frac{x_{2}}{x_{1}}\right)^{\frac{i}{N}}\delta\left(\frac{x_{2}}{x_{1}}\right).
\end{eqnarray}
Then we may multiply both sides by
$\iota_{x_{1};0}(f(x_{1})^{-1})\iota_{x_{2};0}(f(x_{2})^{-1})$ to
get (\ref{ecommutatornew}), noticing that the associativity law
holds in this case.
\end{proof}

The following is the key factorization result:

\bt{tanalyticmodule} Let $\hat{\g}[\sigma]$ be a twisted affine Lie
algebra and let $\pi: \hat{\g}[\sigma]\rightarrow \End W$ be a
integrable representation in the category ${\mathcal{C}}$. Define
linear maps $\pi_{\mathcal{R}}$ and $\pi_{\mathcal{E}}$ from
$\hat{\g}[\sigma]$ to $\End W$ in terms of generating functions by
\begin{eqnarray}
& &\pi_{\mathcal{R}}(a(x)+\alpha {\bf k})
=\psi_{\mathcal{R}}(\pi (a(x)))+\alpha \pi ({\bf k}),\\
& &\pi_{\mathcal{E}}(a(x)+\beta {\bf k})=\psi_{\mathcal{E}}(\pi (a(x)))
\end{eqnarray}
for $a\in \g,\; \alpha,\beta\in \C$, where we extend $\pi$ to
$\hat{\g}[\sigma][[x,x^{-1}]]$ canonically.
Then
\begin{eqnarray}\label{epi-projections}
\pi =\pi_{\mathcal{R}}+\pi_{\mathcal{E}}
\end{eqnarray}
and the linear map
\begin{eqnarray}
\hat{\g}[\sigma]\oplus \hat{\g}[\sigma] \rightarrow \End W:
(u,v)\mapsto \pi_{\mathcal{R}}(u)+\pi_{\E}(v)
\end{eqnarray}
defines a representation of $\hat{\g}[\sigma]\oplus
\hat{\g}[\sigma]$ on $W$. If $(W,\pi)$ is irreducible, we have that
$W$ is an irreducible $\hat{\g}[\sigma]\oplus
\hat{\g}[\sigma]$-module. Furthermore, $(W,\pi_{\mathcal{R}})$ is a
restricted integrable $\hat{\g}[\sigma]$-module (in the category
${\mathcal{R}}$) and $(W,\pi_{\E})$ is a integrable
$\hat{\g}[\sigma]$-module in the category $\E$.\et

\begin{proof} The relation (\ref{epi-projections}) follows from
Proposition \ref{pdecomposition}. It follows immediately from the
defining commutator relations (\ref{edef-commutator}) and Lemma
\ref{lcommutatorrelation} that $(W,\pi_{\mathcal{R}})$ is a
$\hat{\g}[\sigma]$-module and it is clear that it is a restricted
module. (We view ${\bf k}$ as an element of $\bar{\E}(W)$.)
Consequently, $(W,\pi_{\E})$ is a $\hat{\g}[\sigma]$-module, since
$\pi_{\E}=\pi -\pi_{\CR}$.

Let $0\ne f(x)\in \C[x]$ be such that $f(x)\pi(a(x))\in \Hom
(W,W((x^{\frac{1}{N}})))$ for all $a\in \g$. Then
$$f(x)\psi_{\CR}\pi (a(x))=f(x)\pi (a(x)),\;\;\;\;f(x)\psi_{\E}\pi (a(x))=0,$$
so that
\begin{eqnarray}
& &f(x)\pi_{\CR} (a(x))=f(x)\psi_{\CR}\pi (a(x))=f(x)\pi (a(x)),\\
 & &f(x)\pi_{\E}(a(x))=0
\end{eqnarray}
for $a\in \g$. From this we have that $(W,\pi_{\E})$ belongs to the
category $\E$.

For $a\in \g_{i},\,b\in \g_{j}$, using the commutator relations
(\ref{edef-commutator}) and the basic delta-function property we
have
\begin{eqnarray}
& &f(x_{1})[\pi_{\CR}(a(x_{1})),\pi_{\E}(b(x_{2}))]\nonumber\\
&=&f(x_{1})[\pi_{\CR}(a(x_{1})),\pi(b(x_{2}))]
-f(x_{1})[\pi_{\CR}(a(x_{1})),\pi_{\CR}(b(x_{2}))]\nonumber\\
&=&f(x_{1})[\pi (a(x_{1})),\pi (b(x_{2}))]
-f(x_{1})[\pi_{\CR}(a(x_{1})),\pi_{\CR}(b(x_{2}))]\nonumber\\
&=&f(x_{1})\pi
([a,b](x_{2}))x_{1}^{-1}\left(\frac{x_{2}}{x_{1}}\right)^{\frac{i}{N}}\delta\left(\frac{x_{2}}{x_{1}}\right)
+\<a,b\>\pi ({\bf k})f(x_{1})\frac{\partial}{\partial x_{2}}\left[
x_{1}^{-1}\left(\frac{x_{2}}{x_{1}}\right)^{\frac{i}{N}}\delta\left(\frac{x_{2}}{x_{1}}\right)\right]\nonumber\\
& &-f(x_{1})\pi_{\CR}([a,b](x_{2}))
x_{1}^{-1}\left(\frac{x_{2}}{x_{1}}\right)^{\frac{i}{N}}\delta\left(\frac{x_{2}}{x_{1}}\right)
-\<a,b\>\pi_{\CR} ({\bf k})f(x_{1})\frac{\partial}{\partial
x_{2}}\left[
x_{1}^{-1}\left(\frac{x_{2}}{x_{1}}\right)^{\frac{i}{N}}\delta\left(\frac{x_{2}}{x_{1}}\right)\right]
\nonumber\\
&=&f(x_{2})\pi
([a,b](x_{2}))x_{1}^{-1}\left(\frac{x_{2}}{x_{1}}\right)^{\frac{i}{N}}\delta\left(\frac{x_{2}}{x_{1}}\right)
+\<a,b\>\pi ({\bf k})f(x_{1})\frac{\partial}{\partial x_{2}}
\left[x_{1}^{-1}\left(\frac{x_{2}}{x_{1}}\right)^{\frac{i}{N}}\delta\left(\frac{x_{2}}{x_{1}}\right)\nonumber\right]\\
& &-f(x_{2})\pi_{\CR}([a,b](x_{2}))
x_{1}^{-1}\left(\frac{x_{2}}{x_{1}}\right)^{\frac{i}{N}}\delta\left(\frac{x_{2}}{x_{1}}\right)
-\<a,b\>\pi ({\bf k})f(x_{1})\frac{\partial}{\partial x_{2}}
\left[x_{1}^{-1}\left(\frac{x_{2}}{x_{1}}\right)^{\frac{i}{N}}\delta\left(\frac{x_{2}}{x_{1}}\right)\right]
\nonumber\\
&=&0,
\end{eqnarray}
noticing that
$f(x_{2})\pi([a,b](x_{2}))=f(x_{2})\pi_{\CR}([a,b](x_{2})$. Since
$\pi_{\CR}(a(x_{1}))\in \Hom (W,W((x_{1}^{\frac{1}{N}})))$, we can
multiply both sides by $\iota_{x_{1},0}1/f(x_{1})$ and use
associativity to get
\begin{eqnarray}
[\pi_{\CR}(a(x_{1})),\pi_{\E}(b(x_{2}))]=0.
\end{eqnarray}
It follows that $(u,v)\mapsto \pi_{\CR}(u)+\pi_{\E}(v)$ defines a
representation of $\hat{\g}[\sigma]\oplus \hat{\g}[\sigma]$ on $W$.
With $\pi=\pi_{\CR}+\pi_{\E}$, it is clear that if $(W,\pi)$ is
irreducible, $W$ is an irreducible $\hat{\g}[\sigma]\oplus
\hat{\g}[\sigma]$-module.

To prove that $(W,\pi_{\mathcal{R}})$ and $(W,\pi_{\E})$ are
integrable $\hat{\g}$-modules, we have to prove that for $a\in
\g_{\alpha}$ with $\alpha\in \Delta$ and for $n\in \Z$,
$\tilde{a}(n)$ and $\check{a}(n)$ act locally nilpotently on $W$,
where $\pi_{\mathcal{R}}(a\otimes t^{n})=\tilde{a}(n)$ and
$\pi_{\E}(a\otimes t^{n})=\check{a}(n)$.

Let $a\in \g_{\alpha}$ with $\Delta$ and $n\in \Z$. Notice that
$[a(r),a(s)]=0$ for $r,s\in \Z$, since $[a,a]=0$ and $\<a,a\>=0$.
For $w\in W$, we have
$$a(r)\tilde{a}(x)w=a(r)\iota_{x;0}(1/f(x)) (f(x)a(x)w)
=\iota_{x;0}(1/f(x)) (f(x)a(x)a(r)w)=\tilde{a}(x)a(r)w.$$ Thus
\begin{eqnarray}\label{ecommuting-tilde}
a(r)\tilde{a}(s)=\tilde{a}(s)a(r)\;\;\;\mbox{ for }r,s\in \Z.
\end{eqnarray}
Let $w\in W$ be an arbitrarily fixed vector. By Lemma
\ref{lconnection},
$$\tilde{a}(n)w=\sum_{i=0}^{r}\beta_{i}a(n+i)w$$
for some positive integer $r$ and for some complex numbers
$\beta_{1},\dots,\beta_{r}$. Using (\ref{ecommuting-tilde}) we get
\begin{eqnarray}\label{eexp-formula}
\tilde{a}(n)^{p}w=(\beta_{0}a(n)+\cdots +\beta_{r}a(n+r))^{p}w
\;\;\;\mbox{ for any }p\ge 0.
\end{eqnarray}
Since $(W,\pi)$ is an integrable $\hat{\g}$-module, there is a
positive integer $k$ such that
$$a(m)^{k}w=0\;\;\;\mbox{ for }m=n,n+1,\dots,n+r.$$
Combining this with (\ref{eexp-formula}) we obtain
$\tilde{a}(n)^{k(r+1)}w=0$.

Since $\check{a}(n)=a(n)-\tilde{a}(n)$ and $[a(n),\tilde{a}(n)]=0$,
we get
\begin{eqnarray}
\check{a}(n)^{k(r+2)}w=(a(n)-\tilde{a}(n))^{k(r+2)}w =\sum_{i\ge
0}\binom{k(r+2)}{i}(-1)^{i}a(n)^{k(r+2)-i} \tilde{a}(n)^{i}w=0.
\end{eqnarray}
This proves that $\tilde{a}(n)$ and $\check{a}(n)$ act locally
nilpotently on $W$, completing the proof.
\end{proof}

Furthermore, we  have:

\bp{pidentification-1} Let $\hat{\g}[\sigma]$ be a twisted affine
Lie algebra and let $(W_{1},\pi_{1})$ and $(W_{2},\pi_{2})$ be
$\hat{\g}[\sigma]$-modules in the category ${\mathcal{C}}$ and let
$\Phi$ be a $\hat{\g}[\sigma]$-module homomorphism (isomorphism)
from $(W_{1},\pi_{1})$ to $(W_{2},\pi_{2})$. Then $\Phi$ is a
$\hat{\g}[\sigma]$-module homomorphism (isomorphism) {}from
$(W_{1},(\pi_{1})_{\CR})$ to $(W_{2},(\pi_{2})_{\CR})$ and a
$\hat{\g}[\sigma]$-module homomorphism (isomorphism) {}from
$(W_{1},(\pi_{1})_{\E})$ to $(W_{2},(\pi_{2})_{\E})$. \ep

\begin{proof} Let $f(x)$ be a nonzero polynomial such that
for every $a\in \g$,
$$f(x)\pi_{1}(a(x))\in \Hom (W_{1},W_{1}((x^{\frac{1}{N}}))),
\;\; f(x)\pi_{2}(a(x))\in \Hom (W_{2},W_{2}((x^{\frac{1}{N}}))).$$
Then we have
$$f(x)\psi_{\CR}(\pi_{1}(a(x)))=f(x)\pi_{1}(a(x)),\;\;
f(x)\psi_{\CR}(\pi_{2}(a(x)))=f(x)\pi_{2}(a(x)),$$
so that
\begin{eqnarray}
& &f(x)(\pi_{1})_{\CR}(a(x))=
f(x)\psi_{\CR}(\pi_{1}(a(x)))=f(x)\pi_{1}(a(x)),\\
& &f(x)(\pi_{2})_{\CR}(a(x))=
f(x)\psi_{\CR}(\pi_{2}(a(x)))=f(x)\pi_{2}(a(x)).
\end{eqnarray}
For $a\in \g,\; w_{1}\in W_{1}$, we have
\begin{eqnarray}
& &f(x)\Phi((\pi_{1})_{\CR}(a(x))w_{1})=f(x)\Phi(\pi_{1}(a(x))w_{1})
=f(x)\pi_{2}(a(x))\Phi (w_{1})\nonumber\\
& &\hspace{1cm}=f(x)(\pi_{2})_{\CR}(a(x))\Phi (w_{1}).
\end{eqnarray}
Since $\Phi((\pi_{1})_{\CR}(a(x))w_{1}),\; (\pi_{2})_{\CR}(a(x))\Phi
(w_{1})\in W_{2}((x))$, we have
\begin{eqnarray}
\Phi ((\pi_{1})_{\CR}(a(x))w_{1})
=(\pi_{2})_{\CR}(a(x))\Phi (w_{1})\;\;\;\mbox{ for }a\in \g.
\end{eqnarray}
This proves that $\Phi$ is a $\hat{\g}$-module homomorphism {}from
$(W_{1},(\pi_{1})_{\CR})$ to $(W_{2},(\pi_{2})_{\CR})$. Because
$(\pi_{i})_{\E}=\pi_{i}-(\pi_{i})_{\CR}$ for $i=1,2$, it follows
that $\Phi$ is also a $\hat{\g}$-module homomorphism {}from
$(W_{1},(\pi_{1})_{\E})$ to $(W_{2},(\pi_{2})_{\E})$.
\end{proof}

Now, we are in a position to prove our main result:

\bt{tclassification-simple}
Let $\hat{\g}$ be a standard affine Lie algebra.
Every irreducible integrable $\hat{\g}$-module
in the category ${\mathcal{C}}$
is isomorphic to a module of the form
$M\otimes U_{1}(z_{1})\otimes \cdots \otimes U_{r}(z_{r})$,
where $M$ is an irreducible integrable highest weight
$\hat{\g}$-module and $U_{1},\dots,U_{r}$ are finite-dimensional
irreducible $\g$-modules with $z_{1},\dots,z_{r}$ distinct nonzero
complex numbers.
\et

\begin{proof}
Let $\pi: \hat{\g}\rightarrow \End W$ be an irreducible integrable representation of $\hat{\g}$
in the category ${\mathcal{C}}$.
By Theorem \ref{tanalyticmodule}, $W$ is an irreducible
$\hat{\g}\oplus \hat{\g}$-module with $(u,v)$ acting as
$\pi_{\mathcal{R}}(u)+\pi_{\mathcal{E}}(v)$
for $u,v\in \hat{\g}$ and we have $\pi=\pi_{\mathcal{R}}+\pi_{\mathcal{E}}$.
Furthermore, by Proposition \ref{pintegrability},
$(W,\pi_{\mathcal{R}})$ is an integrable restricted $\hat{\g}$-module
and $(W,\pi_{\E})$
is an integrable $\hat{\g}$-module in the category $\E$.
In view of Theorem \ref{tdlm}, $(W,\pi_{\mathcal{R}})$ is a direct sum
of integrable highest weight (irreducible) $\hat{\g}$-modules.
Now it follows immediately from Lemma \ref{ltensor-decomposition-2}
with $A_{1}=A_{2}=U(\hat{\g})$ (which is of countable dimension)
and from Proposition \ref{pevaluation-simple-modules}.
\end{proof}


\begin{thebibliography}{FKRW}

\bibitem[B]{b}
R. E. Borcherds, Vertex algebras, Kac-Moody algebras, and the Monster,
{\it Proc. Natl. Acad. Sci. USA} {\bf 83} (1986), 3068-3071.

\bibitem[C]{ch}
V. Chari, Integrable representations of affine Lie algebras,
{\em Invent. Math.} {\bf 85} (1986), 317-335.

\bibitem[CP1]{cp1}
V. Chari and A. N. Pressley, New unitary representations of loop groups,
{\em Math. Ann.} {\bf 275} (1986), 87-104.

\bibitem[CP2]{cp2}
V. Chari and A. N. Pressley, A new family of irreducible, integrable modules
for affine Lie algebras, {\em Math. Ann.} {\bf 277} (1987), 543-562.

\bibitem[CP3]{cp3}
V. Chari and A. N. Pressley, Integrable representations of twisted
affine Lie algebras, {\em J. Algebra} {\bf 113} (1988), 438-464.

\bibitem[FLM]{flm}
I. Frenkel, J. Lepowsky and A. Meurman, {\it Vertex Operator
Algebras and the Monster}, Pure and Appl. Math., {\bf Vol. 134},
Academic Press, Boston, 1988.

\bibitem[H]{hum}
J. Humphreys, {\em Introduction to Lie Algebras and Representation
Theory,} Springer-Verlag, New York, 1972.

\bibitem[K]{kacbook}
V. G. Kac, Infinite Dimensional Lie Algebras,
3rd edition, Cambridge University Press, 1990.

\bibitem[Li]{li}
H.-S. Li, On certain categories of modules for affine Lie algebras,
{\em Math. Z.} {\bf 248} (2004), 635-664.

\bibitem[LL]{ll}
J. Lepowsky and H.-S. Li,
{\em Introduction to Vertex Operator Algebras and Their Representations},
Progress in Math., {\bf Vol.} 227, Birkh\"auser, Boston, 2004.



\end{thebibliography}
\end{document}